\documentclass[12pt] {amsart}
\usepackage{comment}
\usepackage{enumerate}
\usepackage{amssymb, amsmath}

\newtheorem{thm}{Theorem}[section]
\newtheorem{prop}[thm]{Proposition}
\newtheorem{lemma}[thm]{Lemma}

\usepackage{url}

\def\Z{{\mathbb Z}}

\usepackage{hyperref}
\usepackage{geometry}
\usepackage[T1]{fontenc}
\usepackage{hyperref}
\usepackage{longtable}
\usepackage{enumerate}

 \hypersetup{
   linkcolor=green,
    filecolor=black,      
    urlcolor= cyan,
    }

\bibliography{PD.bib}

\begin{document}

\title{Quartic Index Form Equations and Monogenizations of Quartic Orders }

\author{Shabnam Akhtari}

 \email {akhtari@uoregon.edu}

\subjclass[2000]{11D25, 11D45,  11R04, 11R16}
\begin{abstract}
 Some upper bounds for the number of monogenizations of quartic orders are established by considering certain classical Diophantine equations, namely index form equations in quartic number fields, and cubic and quartic Thue equations.
\end{abstract}
\keywords{Monogenizations of a quartic order, Index form equations, Thue equations}

\maketitle

\section{Introduction}

Let $K$ be an algebraic number field and $O_K$ its ring of integers. Let $O$ be an order in  $K$ (a subring of $O_K$ with quotient field $K$). We call the ring $O$ {\it monogenic} if it is generated by one element as a $\Z$-algebra, i.e., $O=\Z[\alpha]$ for some $\alpha\in O$; the element $\alpha$ is  called a {\it monogenizer} of $O$.  If $\alpha$ is a monogenizer of $O$, than so is $\pm \alpha+c$ for any $c\in\Z$. We call two monogenizers $\alpha$ and $\alpha'$ of $O$  {\it equivalent} if $\alpha'=\pm \alpha +c$ for some $c\in\Z$.  Then by a
 {\it monogenization} of $O$, we mean
an equivalence class of monogenizers of $O$. By fundamental work of Gy\H{o}ry \cite{Gyo76}, we know that any order in an algebraic number field can have at most finitely many monogenizations and that  effectively computable upper bounds on the number of these monogenizations can be determined.  It is a difficult computational problem to find or even count the monogenizations of a given order (many computational examples, interesting special cases and  efficient  algorithms in low degree number fields may be found in \cite{thebook}).

We are interested in counting the number of  monogenizations of a given order.
An overview of various results on estimates for the number of monogenizations  of orders in number fields is given in  \cite{Eve11}. There are further extensions and generalizations of such results in  \cite{EGbook}  (in particular, see the relevant results  in Section 9.1).

Monogenicity of algebraic number rings has a long history. It is an interesting problem to decide whether a given number field $K$ is monogenic, that is, whether its ring of integers $O_K$, which is the maximal order in $K$,  is monogenic. It is well known that quadratic number fields are monogenic. In 1878 Dedekind \cite{Ded} gave the first example of a non-monogenic cubic field. It is an open conjecture that most of number fields of degree greater than $2$ are not monogenic. For recent progress in this direction in the cases of cubic and quartic number fields, we refer the reader to the work of Alp\"{o}ge, Bhargava,  and Shnidman in \cite{ABS1, ABS2}.

In this article we focus on the problem of counting the number of monogenizations of a quartic order. In \cite{EG85} Evertse and Gy\H{o}ry proved explicit upper bounds for the number of monogenizations of an order in a number field $K$. These bounds  depend only on the degree of $K$.
 The best known result for $n\geq 4$ is due to Evertse, who proved in \cite{Eve11} that an order $O$ in a number field $K$ of degree $n$ can have at most
 $2^{4(n+5)(n-2)}$
monogenizations. In the case $n=4$ , Evertse’s result  shows that an order in a quartic field can have at most $2^{72}$ monogenizations.  Recently Bhargava  gave an improved bound  in \cite{Bha-notes}, showing that an order in a quartic number field can have at most $2760$ monogenizations (and even fewer when the discriminant of the order is large enough). We give another proof for Theorem 1.1 of  \cite{Bha-notes}. 

\begin{thm}\label{whatisM}
 Let $O$ be an order in a quartic number field. The number of monogenizations of  $O$  is at most $2760$. If the absolute value of the  discriminant of $O$ is sufficiently large,  the number of monogenizations of  $O$  is at most  $182$. Moreover, if the discriminant of $O$ is negative and has sufficiently large absolute value, the number of monogenizations of  $O$  is at most  $70$.
\end{thm}

In the above theorem the assumptions about the size of the discriminant are the result of such assumptions to overcome certain technical difficulties in some approximation methods used to prove Propositions \ref{BenOka}, \ref{main26}, and  \ref{quartic-negative}. These restrictions can be expressed explicitly. For instance,  assuming the absolute value of the discriminant is at least $10^{500}$ will suffice (see \cite{Akh10, A}, where such explicit values are established but no effort has been made to optimize them).  It is known  that there are only finitely many quartic number fields with the absolute value of their discriminants bounded by a constant (see \cite{BiMe,EG91}). By the identity in \eqref{discalphaiK}, which relates the discriminant of an order to that of the underlying number field, Theorem \ref{whatisM} implies that  with at most finitely many exceptions, a quartic order with positive discriminant can have at most $182$ monogenizations and a quartic order with negative discriminant can have at most $70$ monogenizations.

Our approach involves refining and modifying an algorithmic  method developed in 1996 by  Ga\'{a}l,  Peth\H{o} and Pohst in \cite{GaPePo} to solve an index form equation $I(X, Y, Z) = \pm 1$ in a quartic number field. Using this method, we will be able to associate explicit polynomials and binary  and ternary forms to a monogenic order and a fixed monogenizer of that, and eventually reduce our problem to the resolution of a number of Thue equations of degree $3$ and $4$. The proof in \cite{Bha-notes} uses a more abstract viewpoint by utilizing  two ways of  parametrizing quartic rings, one established by Bhargava in \cite{Bha04} and another one established by Wood in \cite{Woo12}.

\section{Preliminaries:  Discriminants, Thue Equations, Discriminant  and Index Form Equations}

\subsection{Discriminants}

We recall the definitions of  discriminants of orders, polynomials, algebraic numbers,  and binary forms which will be frequently used throughout this manuscript. We will also refer to the discriminant of number fields. The discriminant of a number field $K$ is the discriminant of its maximal order, the ring of integers $O_K$. For
 $K =\mathbb{Q}(\alpha)$, the discriminant of $K$ can be expressed in terms of the discriminant of  the algebraic number $\alpha$ and its index in $\mathbb{Q}(\alpha)$. The index of an algebraic integer and the discriminant of orders are defined in \textsection  \ref{DisIn}.

Let $\mathbf{P}(T) \in \mathbb{Z}[T]$ be a polynomial of degree $n$ and leading coefficient
 $a \in \mathbb{Z}$.
The discriminant $\textrm{Disc}(\mathbf{P})$ of $\mathbf{P}(T)$ is  
$$
\textrm{Disc}(\mathbf{P}) = a^{2n -2} \prod_{i <j} (\gamma_i - \gamma_j)^2,
$$
where $\gamma_1, \ldots, \gamma_n \in \mathbb{C}$ are the roots of $\mathbf{P}(T)$.

The discriminant of an algebraic number is defined as the discriminant of its minimal polynomial.

Let $F(U, V) \in \mathbb{Z}[U, V]$ be a binary form of degree $n$ that factors over $\mathbb{C}$ as
$$
\prod_{i=1}^{n} (\alpha_{i} U - \beta_{i} V).
$$
The discriminant $D(F)$ of $F$ is given by
\begin{equation}\label{defofdisc}
D(F) = \prod_{i<j} (\alpha_{i} \beta_{j} - \alpha_{j}\beta_{i})^2.
\end{equation}

We note that the discriminant of the polynomial $F(U, 1) \in \mathbb{Z}[U]$ is equal to the discriminant of the binary form $F(U, V)  \in \mathbb{Z}[U, V]$.

\subsection{Discriminant and Index Form Equations}\label{DisIn}

Let $K$ be an algebraic number field of degree $n$. 
Let $\alpha_{1}, \ldots, \alpha_{n}$  a linearly independent set of $n$ elements of $K$. Let $\sigma_{1}, \ldots, \sigma_{n} : K \rightarrow \mathbb{C}$ be all the
embeddings of $K$ into $\mathbb{C}$. The discriminant of $(\alpha_{1}, \ldots, \alpha_{n})$ is defined as the square of the determinant of an $n \times n$ matrix;
 $$
D_{K/\mathbb{Q}}(\alpha_{1}, \ldots, \alpha_{n}): = \left( \textrm{det}(\sigma_{i}(\alpha_{j})\right)^2,
$$
where $i, j \in \{1, \ldots, n\}$.

 If $\{\beta_{1}, \ldots, \beta_{n}\}$ forms a basis for $O_K$, then the discriminant of $K$ is
 $$
 D_{K} =  D_{K/\mathbb{Q}}(\beta_{1}, \ldots, \beta_{n}).
 $$

 Let $\gamma_1, \gamma_2, \ldots, \gamma_n$ be an integral basis for an order $O$ in a number field  $K$ of degree $n$ (we note that by definition an order is a full-rank  $\mathbb{Z}$-module in $O_K$). The discriminant of $O$ is defined as
 $
D_{K/\mathbb{Q}}(\gamma_{1}, \ldots, \gamma_{n})
$
 and is independent of the choice of the integral basis  $\gamma_1, \gamma_2, \ldots, \gamma_n$ (see \cite{Koch}, or any  introductory text in algebraic number theory).

The following basic well-known  lemmas are due to Hensel and can be found in \cite{Hen}.

\begin{lemma}\label{NarLemma}
Let $\alpha_{1}, \ldots, \alpha_{n} \in {O}_{K}$ be linearly independent over $\mathbb{Q}$  and set 
$$
{O} = \mathbb{Z}[\alpha_{1}, \ldots, \alpha_{n}]. 
$$
then 
\begin{equation}\label{discalphaiK}
D_{K/\mathbb{Q}}(\alpha_{1}, \ldots, \alpha_{n}) = J^2 D_{K},
\end{equation}
where 
${O}_{K}^{+}$ and ${O}^{+}$ are the additive groups of the modules ${O}_{K}$ and  ${O}$, respectively, and $
J = ({O}_{K}^{+} : {O}^{+})
$ is the module index.
\end{lemma}

For every $\gamma \in K$,  we denote the algebraic  conjugates of $\gamma$ by  $\gamma^{(i)}$  ($1\leq i \leq n$). Let $\{1, \omega_{2}, \ldots, \omega_{n}\}$ be an integral basis of $K$. Let
$$
\mathbf{X} = (X_{1} ,  \ldots , X_{n}),
$$
and
\begin{equation}\label{Ldef}
L(\mathbf{X}) = X_{1} + \omega_{2} X_{2} + \ldots +  \omega_{n} X_{n},
\end{equation}
with algebraic conjugates
$$
L^{(i)} (\mathbf{X}) = X_{1} + \omega_{2}^{(i)} X_{2} + \ldots +  \omega_{n}^{(i)} X_{n},
$$
($1\leq i \leq n$).
Kronecker and Hensel  called the form $L(\mathbf{X})$ the \emph{Fundamentalform} and 
\begin{equation}\label{Ddef}
D_{K/\mathbb{Q}} (L(\mathbf{X}))=  \prod_{1\leq i < j \leq n}\left(  L^{(i)} (\mathbf{X})  - L^{(j)} (\mathbf{X})   \right)^2
\end{equation}
the \emph{Fundamentaldiskriminante}.
\begin{lemma}\label{DIDK}
We have
$$
D_{K/\mathbb{Q}} (L(\mathbf{X}))=   \left( I(X_{1} ,  \ldots , X_{n})  \right)^2 D_{K},
$$
where $D_{K}$ is the discriminant of the field $K$, the linear form  $L(\mathbf{X})$ and its discriminant  are defined in \eqref{Ldef} and \eqref{Ddef}, and $I(X_{1} ,  \ldots , X_{n})  $ is a homogeneous form in $n-1$ variables of degree $\frac{n(n-1)}{2}$ with integer coefficients. 
\end{lemma}

The form $I(X_{1} ,  \ldots , X_{n})$ in the statement of Lemma \ref{DIDK}  is called the index form corresponding to the integral basis $\{1, \omega_{2}, \ldots, \omega_{n}\}$.
An important property of the index form is that for any  algebraic integer  
$$
\alpha = x_{1} + x_{2} \omega_{2} + \ldots + x_{n} \omega_{n},
$$
with $K = \mathbb{Q}(\alpha)$, by Lemma \ref{DIDK}
we have
$$
I(\alpha) = \left| I(x_{2}, \ldots, x_{n}) \right|,
$$
where $I(\alpha)$ is the  index of the module $\mathbb{Z}[\alpha]$ in $O_K$.
The index form is independent of the variable $X_{1}$, for if $\beta = \alpha + a$, where $a \in \mathbb{Z}$, then $I(\alpha) = I(\beta)$.

We remark that in a  cubic number field an index form equation is in fact a cubic Thue equation (see \textsection \ref{ThueBounds} for the definition)
\begin{equation*}\label{mainindexformeq3}
I(X_{2}, X_{3}) = \pm m,
\end{equation*}
where $m \in \mathbb{Z}$.  
In \cite{Akh-Comb} we have discussed some results about cubic Thue equations and their consequences in resolving index form equations and counting the number of monogenizations of a cubic ring.

\subsection{Upper Bounds on the Number of Solutions of Cubic and Quartic Thue Equations}\label{ThueBounds}

Let $F(U , V) \in \mathbb{Z}[U , V]$ be a binary form of degree at least $3$. If $F(U, V)$ is irreducible over $\mathbb{Q}$, for any integer $m$, it is shown in \cite{Thu}  that the equation 
$$
F(U , V) = m
$$
has at most finitely many solutions in integers $U$, $V$. These equations are called \emph{Thue equations}.  We will summarize some useful results on the number of integer solutions of binary cubic and quartic Thue equations. In Propositions  \ref{BenOka},  \ref{Del-Nag}, \ref{main26}, and \ref{quartic-negative},  two pairs of solutions 
$(u, v), (-u, -v) \in \mathbb{Z}^2$ are considered as one solution.

The following is the combination of main results due to Bennett in \cite{B} and Okazaki in \cite{Oka} (see also \cite{A}).
\begin{prop}\label{BenOka}
 A cubic Thue equation $F(U, V) = \pm 1$ has at most $10$ integer solutions.
  If the absolute value of the discriminant of $F(U , V)$ is sufficiently  large then $F(U, V) = \pm 1$ has at most $7$ integer solutions.
\end{prop}

The following result is established independently by Delone and Nagell in \cite{Del} and \cite{Nag}, respectively. 
\begin{prop}\label{Del-Nag}
 Let $F(U, V) \in \mathbb{Z}[U, V]$ be a cubic binary form with negative discriminant. The Thue equation $F(U, V) = \pm 1$ has at most $5$  integers solutions. 
\end{prop}

The following is Theorem  A.1 of
\cite{Bha-notes}, where   results from \cite{Akhsmallm, Akh10, Mike}  are combined to obtain upper bounds for the number of integral solutions to quartic Thue equations.
\begin{prop}\label{main26}
A quartic Thue equation $F(U , V) = \pm 1$ has at most $276$  integer solutions. If the absolute value of the discriminant of $F(U , V)$ is sufficiently  large then the quartic Thue equation $F(U , V) = \pm 1$ has at most $26$  integer solutions.
\end{prop}

The following is part of the main theorem in \cite{Akh10}.
\begin{prop}\label{quartic-negative}
Let $F(U, V) \in \mathbb{Z}[U, V]$ be a quartic binary form with negative discriminant. If the absolute value of the  discriminant of $F(U , V)$ is sufficiently  large, the Thue equation $F(U , V) = \pm 1$ has at most $14$  integer solutions.
\end{prop}

\subsection{Matrix Actions on Binary Forms}

We summarize some trivial facts about matrix actions on binary forms that are well known to those in the field. 
 Let $F(U, V) \in \mathbb{Z}[U, V]$ and
$
A = \left( \begin{array}{cc}
a & b \\
c & d \end{array} \right)
$ be a $2 \times 2$ matrix with integer entries. 
We define the binary form $$F_{A}(U, V) \in \mathbb{Z}[U, V]$$
 by
$$
F_{A}(U , V) = F(a U+ b V,  c U+ d V).
$$

Via the definition  \eqref{defofdisc}, we 
observe that for any $2 \times 2$ matrix $A$ with integer entries
\begin{equation}\label{St6}
D(F_{A}) = (\textrm{det} A)^{n (n-1)} D(F).
\end{equation}

We say that two integral binary forms $F$ and $G$ are {\it equivalent} if $G = \pm F_{A}$ for some $A \in \textrm{GL}_{2}(\mathbb{Z})$. This is in fact an equivalence relationship. Moreover, the discriminants of two equivalent forms are equal.

For  $A = \left( \begin{array}{cc}
a & b \\
c & d \end{array} \right)
 \in \textrm{GL}_{2}(\mathbb{Z})$, and any $(u, v) \in \mathbb{Z}^2$, we clearly have 
$$A^{-1} = \pm \left( \begin{array}{cc}
d& -b \\
-c & a \end{array} \right)
$$
and 
$$F_{A} ( du-bv , -cu+av) =  \pm F(u , v).
$$
Therefore, there is a one-to-one correspondence between the possible solutions of the Thue equation $F(U , V) = \pm 1$ and those of the Thue equation $F_A(U , V) = \pm 1$.


\section{Index Form Equations in Quartic Number Fields}\label{GPP-section}

Let $\xi$ be a quartic algebraic integer with the minimal polynomial 
\begin{equation}\label{minpolydef}
\mathbf{P}(T) = T^4 + a_{1} T^3 + a_{2} T^2 + a_{3} T + a_{4} \in \mathbb{Z}[T].
\end{equation}
Let  $K = \mathbb{Q}(\xi)$.
Suppose that  $\omega_{1} = 1$, $\omega_{2}$, $\omega_{3}$ and $\omega_{4}$ form an integral basis for the quartic number field $K$.  We write $\sigma_1$, $\sigma_2$, $\sigma_3$ and $\sigma_4$ for the distinct  embeddings of $K$ into $\mathbb{C}$.
For $i =1, 2, 3, 4$, we define the linear forms
$$
 l_{i}(X, Y, Z) = X \omega^{(i)}_{2} + Y \omega^{(i)}_{3} + Z \omega^{(i)}_{4},
$$
where $\omega^{(i)}_{j} =  \sigma_i\left(\omega_{j}\right)$.

The \emph{discriminant form} corresponding to the  integral basis  $\{1, \omega_{2}, \omega_{3}, \omega_{4}\}$ is defined by
$$
D_{K/\mathbb{Q}}(X \omega_{2} + Y \omega_{3} + Z \omega_{4}) = \prod_{1\leq i < j \leq 4} \left(l_{i}(X , Y, Z) - l_{j}(X , Y, Z)  \right)^2.
$$
We have
\begin{equation}\label{DID}
D_{K/\mathbb{Q}}(X \omega_{2} + Y \omega_{3} + Z \omega_{4}) = \left( I(X, Y, Z) \right)^2 D_K, 
\end{equation}
where $D_K$ is the discriminant of the number field $K$ and $I(X, Y, Z) \in \mathbb{Z}[X, Y, Z]$ is the \emph{index form} corresponding to the fixed integral basis  $\{1, \omega_{2}, \omega_{3}, \omega_{4}\}$.  The integral  ternary form  $I(X, Y, Z)$ has degree $6$. For any algebraic integer $\alpha = a + x\omega_{2} + y \omega_{3} +z \omega_{4}$, with $a, x, y, z \in \mathbb{Z}$, the index $I(\alpha)$ is equal to
 $\left| I(x , y, z) \right|$, where $I(\alpha)$ is the module index of $\mathbb{Z}[\alpha]$ in $O_{K}$, the ring of integers of $K$.
In this section we consider the \emph{index form equation} 
\begin{equation}\label{mainindexformeq}
I(X, Y, Z) = \pm m
\end{equation}
where $m \in \mathbb{Z}$.

We follow a simple and efficient algorithm given by Ga\'al, Peth\H{o} and Pohst in 
  \cite{GaPePo}, where they reduce the problem of solving an index form equation  in  a quartic number field to the problem of finding all solutions $(u_i, v_i) \in \mathbb{Z}^2$ of a cubic Thue equation $F(U , V) = \pm h$, with $h \in \mathbb{Z}$, and  the resolution of corresponding systems of quadratic equations $\mathbf{Q}_{1}(X , Y, Z) = u_i$, $\mathbf{Q}_{2}(X , Y,  Z) = v_i$, where $F(U, V) \in \mathbb{Z}[U, V]$  is a cubic form, and $\mathbf{Q}_{1}(X, Y, Z)$ and $\mathbf{Q}_{2}(X, Y, Z)$ are integral ternary quadratic forms.
We state this reduction more precisely in Proposition \ref{GPP}.

We denote by $I_{0}$ the index of the algebraic integer  $\xi$. Then 
$$
I_{0} = I(\xi) = \left|I(x_0, y_0, z_0)\right|,
$$ 
where $\xi = a_\xi + x_0 \omega_{2} + y_0 \omega_{3} + z_0 \omega_{4}$, and $a_\xi, x_0, y_0, z_0 \in \mathbb{Z}$. Once again we remark that the algebraic integers $\xi$ and $\xi- a_\xi$ have the same index in $O_K$. Since $I_0$ is the index of $\mathbb{Z}[\xi]$ in $O_{K}$, for every algebraic integer $\alpha \in O_{K}$, we have
$$
I_0 \alpha \in \mathbb{Z}[\xi].
$$

Assume that $(x_{1}, y_1, z_1) \in \mathbb{Z}^3$ satisfies \eqref{mainindexformeq}. Let 
\begin{equation}\label{alphainomega}
\alpha = x_{1} \omega_2+ y_1 \omega_{3} + z_1 \omega_{4},
\end{equation}
and
\begin{equation}\label{I0xyz}
\alpha' = I_0 \alpha = a'_{\alpha} + x'_1 \xi + y'_{1}\xi^2 + z'_{1} \xi^3 \in \mathbb{Z}[\xi].
\end{equation}
We have 
\begin{equation}\label{I0indexformeq}
I(\alpha') = I(x' \xi + y' \xi^2 + z' \xi^3) = \pm I_{0}^6 m.
\end{equation}

We denote by $\xi^{(i)}$ and $\alpha'^{(i)}$ the  algebraic conjugates of $\xi$ and $\alpha'$, for $i = 1, 2, 3, 4$. 
 Dividing both sides of the equation \eqref{I0indexformeq} by $I(\xi) = I_{0}$, we obtain
\begin{equation}\label{prod1}
\prod_{(i, j, k, l)} \left(\frac{\alpha'^{(i)} - \alpha'^{(j)}}{\xi^{(i)} - \xi^{(j)}}\right) \left(\frac{\alpha'^{(k)} - \alpha'^{(l)}}{\xi^{(k)} - \xi^{(l)}}\right) = \pm \frac{I_{0}^6 m}{I_{0}} = \pm I_{0}^5 m,
\end{equation}
where the above product is taken for $(i, j, k, l) = (1, 2, 3, 4), (1, 3, 2, 4), (1, 4, 2, 3)$. For each $(i, j, k, l)$, via \eqref{I0xyz}, we have
\begin{equation}\label{ijkl}
 \left(\frac{\alpha'^{(i)} - \alpha'^{(j)}}{\xi^{(i)} - \xi^{(j)}}\right) \left(\frac{\alpha'^{(k)} - \alpha'^{(l)}}{\xi^{(k)} - \xi^{(l)}}\right) = \mathbf{Q}_{1}(x'_1 , y'_1, z'_1) - \xi_{i, j, k,l} \mathbf{Q}_{2}(x'_1 , y'_1, z'_1),
\end{equation}
where
$$
\xi_{i, j, k,l} = \xi^{(i)} \xi^{(j)} + \xi^{(k)} \xi^{(l)},
$$
\begin{eqnarray}\label{defofQ1}
&&\, \, \, \, \, \, \,  \mathbf{Q}_{1}(X , Y, Z) =  \\ \nonumber
&&X^2 - a_{1} XY + a_{2} Y^2 + (a_{1}^2 - 2 a_{2}) XZ 
+(a_{3} - a_{1} a_{2}) YZ + (- a_{1} a_{3} + a_{2}^2 + a_{4}) Z^2,
\end{eqnarray}
and
\begin{equation}\label{defofQ2}
\mathbf{Q}_{2}(X , Y, Z) = Y^2 - XZ - a_{1} YZ + a_{2} Z^2.
\end{equation}
The coefficients of the quadratic forms $\mathbf{Q}_{1}(X, Y,Z)$ and $\mathbf{Q}_{2}(X, Y, Z)$ are expressed in terms of the coefficients of $\mathbf{P}(T)$, the minimal polynomial of $\xi$ given in \eqref{minpolydef}.
For  each $(i, j, k, l) = (1, 2, 3, 4), (1, 3, 2, 4), (1, 4, 2, 3)$, we define the linear  form
$$
 \mathcal{P}(i, j, k, l) =   \mathcal{P}(i, j, k, l) (U, V) = U  - \xi_{1, 2, 3, 4} V.
$$
Taking $U = \mathbf{Q}_1(X, Y, Z)$ and  $V = \mathbf{Q}_2(X, Y, Z)$,  by  \eqref{prod1} and \eqref{ijkl}, we obtain
\begin{equation}\label{prod2}
\prod_{(i, j, k, l)}  \mathcal{P}(i, j, k, l) = (U - \xi_{1, 2, 3, 4} V) (U - \xi_{1, 3, 2, 4} V) (U - \xi_{1, 4, 2, 3} V) =  \pm I_{0}^5 m,
\end{equation}
where
the product is taken over $(i, j, k, l) = (1, 2, 3, 4), (1, 3, 2, 4), (1, 4, 2, 3)$.

The left-hand side of \eqref{prod2} is a cubic binary form in $U$ and $V$ whose coefficients are symmetric polynomials of $\xi^{(1)}$, $\xi^{(2)}$, $\xi^{(3)}$, $\xi^{(4)}$. Simple and routine  calculations show that this integral cubic binary form is 
\begin{eqnarray}\label{defofFu,v}
&& \prod_{(i, j, k, l)}  \mathcal{P}(i, j, k, l)(U , V) = F(U, V)\\ \nonumber
& & = U^3 - a_{2} U^2 V+ (a_{1} a_{3} - 4a_{4}) U V^2 + (4 a_{2} a_{4} - a_{3}^2 - a_{1}^2 a_{4}) V^3.
\end{eqnarray} 
 The cubic polynomial 
 $$F(T, 1) = T^3 - a_{2} T^2 + (a_{1} a_{3} - 4a_{4}) T+ (4 a_{2} a_{4} - a_{3}^2 - a_{1}^2 a_{4}) 
 $$
  is called the \emph{cubic resolvent polynomial}  of  $\mathbf{P}(T)$, the minimal polynomial of $\xi$. The discriminant of $\mathbf{P}(T) \in \mathbb{Z}[T]$ is equal to the discriminant of $F(T , 1) \in \mathbb{Z}[T]$ and therefore to the discriminant of $F(U, V) \in \mathbb{Z}[U, V]$. Since the discriminant of the minimal polynomial $\mathbf{P}(T)$ is not zero, we conclude that $F(U , V)$ will factor into three pairwise non-proportional linear factors over $\mathbb{C}$. This, together with \eqref{prod2}, implies that the three cubic algebraic  integers $\xi_{1, 2, 3, 4}$,  $\xi_{1, 3, 2, 4}$, and  $\xi_{1, 4, 2, 3}$ are  distinct   algebraic conjugates over $\mathbb{Q}$.
The above argument can be found in  \cite{thebook} and \cite{GaPePo}, and implies the following.
\begin{prop}\label{GPP}
Let  $\xi$ be a quartic algebraic integer and
$$I_{0} = I(\xi).$$
 Assume that 
 $I(X, Y, Z) \in \mathbb{Z}[X, Y, Z]$ is an index form in the quartic number field $\mathbb{Q}(\xi)$.
The triple $(x , y, z) \in \mathbb{Z}^3$ is a solution of the index form equation 
$$I(X, Y, Z) = \pm m, $$
 with $m \in \mathbb{Z}$,  if and only if there exists a solution $(u, v) \in \mathbb{Z}^2$ of the cubic Thue equation 
\begin{equation}\label{cubicThueeq}
F(U, V) = \pm I_{0}^5 m
\end{equation}
such that $(x , y, z)$ satisfies the system of quadratic ternary equations
\begin{equation}\label{twoquadratics}
\mathbf{Q}_{1}(X , Y, Z) = u , \, \, 
\mathbf{Q}_{2}(X , Y, Z) = v,
\end{equation}
where $F(U, V)$ is an integral cubic binary form and $\mathbf{Q}_{1}(X , Y, Z)$  and $\mathbf{Q}_{2}(X , Y, Z)$ are integral quadratic ternary forms, respectively defined in \eqref{defofFu,v},  \eqref{defofQ1} and \eqref{defofQ2} with coefficients  expressed in terms of the coefficients of the minimal polynomial of the fixed generator $\xi$.
\end{prop}

Proposition \ref{GPP} provides a general algorithm to find algebraic integers with index $m$ in the quartic number field $K$ by fixing any algebraic integer $\xi$ that generates  $K$. So in general the quantities $I_0$ and $m$ need not to be related. Using an argument in Mordell's  work \cite{MorBook},  in  \cite{GaPePo}  it is shown that all solutions of an index form equation in a quartic number field can be found through solving finitely many cubic and quartic Thue equations (see Theorems 1 and 2, as well as equations (8), (9) and (10) of  \cite{GaPePo}). In \textsection \ref{proofsofmono} we will modify the argument in  \cite{GaPePo}  and apply our modification to an index form equation of the shape $I(X, Y, Z) = \pm 1$ connected to the quartic ring generated by an algebraic integer $\xi$. This will enable us to count  these Thue equations more efficiently.  Moreover, it turns out that in this case the right-hand sides of our Thue equations are $\pm1$, and therefore we may  apply absolute upper bounds for the number of integer solutions  recorded in \textsection \ref{ThueBounds}.  This way we can provide an absolute upper bound for the number of solutions of the index form equation that we study. These solutions will correspond to different monogenizations of the quartic order $\mathbb{Z}[\xi]$.

We end this section by recording another important relation between the ternary quadratic forms $\mathbf{Q}_1$ and $\mathbf{Q}_2$, defined in \eqref{defofQ1} and \eqref{defofQ2},  and the integer values represented by the cubic form $F(U, V)$.
For $(u_0, v_0) \in \mathbb{Z}^2$,
we define 
\begin{equation}\label{defofQ0}
\mathbf{Q}(X, Y, Z)= u_0 \mathbf{Q}_{2}(X , Y, Z) - v_0 \mathbf{Q}_{1}(X , Y, Z).
\end{equation}
Let $M_{\mathbf{Q}}$ be the $3 \times 3$ symmetric Gram matrix of the quadratic form $\mathbf{Q}(X, Y, Z)$. We have
\begin{equation}\label{DetQ0=F}
4\left| \textrm{Det} (M_{\mathbf{Q}})\right| = \left|F(u_0 , v_0)\right|,
\end{equation}
where  $F(U, V)$ is defined in \eqref{defofFu,v}.
The identity \eqref{DetQ0=F} can be verified easily and is established as an implication of Lemma 1 of  \cite{GaPePo}.  Its proof can also be found in Lemma 6.1.1 of \cite{thebook}.
The identity \eqref{DetQ0=F} is not used in our proofs, but it is crucial in confirming that the ternary quadratic forms $\mathbf{Q}_1$ and $\mathbf{Q}_2$ form a pair that parametrizes a quartic ring in the sense of Bhargava's work \cite{Bha04}. Such a parametrization is used in Bhargava's proof of Theorem \ref{whatisM}  in \cite{Bha-notes}. Another ingredient in \cite{Bha-notes} is a beautiful parametrization due to Wood in \cite{Woo12} for quartic rings. We do not use any of these two parametrizations. However, in light of identities  \eqref{defofQ0}  and  \eqref{DetQ0=F},   one could view  our discussion in the following section  as an explicit way of expressing polynomials and binary forms that are appearing (implicitly) in  Bhargava's and Wood's methods of parametrization.

\section{Proof of Theorem \ref{whatisM}}\label{proofsofmono}

When treating a general index form equation in a quartic number field, one needs to  consider the identity \eqref{I0xyz} in order to have integer values for $x'_1, y'_1$ and $z'_1$. 
In Theorem \ref{whatisM} we
 are  interested in finding other possible monogenizers for a monogenized ring $\mathbb{Z}[\xi]$. Therefore,  we are  looking for algebraic integers $\alpha \in \mathbb{Z}[\xi]$ that satisfy the index form equation \eqref{mainindexformeq}. In this case, under the assumption  $\alpha \in  \mathbb{Z}[\xi]$, we may express \eqref{I0xyz} as
\begin{equation}\label{noI0xyz}
 \alpha = a_{\alpha} + x \xi + y\xi^2 + z \xi^3,
\end{equation}
with $x, y, z \in \mathbb{Z}$.  This will simplify some of the equations introduced in \textsection  \ref{GPP-section}. Another simple observation is that
if $\mathbb{Z}[\xi] = \mathbb{Z}[\alpha]$, then the algebraic integers $\alpha$ and $\xi$ have the same index in the ring of integers of the underlying number field $\mathbb{Q}(\alpha) = \mathbb{Q}(\xi)$, and therefore  in the index form equations \eqref{I0indexformeq} and \eqref{prod1}, we may take $I_0 = m$.

Let $K$ be a quartic number field and $\xi$  an algebraic integer in $K$ of index $I_0 = m$. We are interested in finding other monogenizers of $\mathbb{Z}[\xi]$. 
After replacing \eqref{I0xyz} by \eqref{noI0xyz},  for $\alpha \in \mathbb{Z}[\xi]$ the identity \eqref{prod1}  becomes
\begin{equation}\label{prod1=1}
\prod_{(i, j, k, l)} \left(\frac{\alpha^{(i)} - \alpha^{(j)}}{\xi^{(i)} - \xi^{(j)}}\right) \left(\frac{\alpha^{(k)} - \alpha^{(l)}}{\xi^{(k)} - \xi^{(l)}}\right) = \pm 1.
\end{equation}
Therefore, in Proposition \ref{GPP}, we may consider the cubic Thue equation 
\begin{equation}\label{cubicthue=1}
F(U, V) = \pm 1.
\end{equation}

In fact, we obtain the following modification of Proposition \ref{GPP}.
\begin{lemma}\label{GPPmodified}
The algebraic integer $x\xi + y \xi^2 + z\xi^3$, with $x , y, z \in \mathbb{Z}$ is a monogenizer of $\mathbb{Z}[\xi]$ if and only if there is a solution $(u, v) \in \mathbb{Z}^2$ of the cubic Thue equation 
\begin{equation}\label{cubicThueeqmodified}
F(U, V) = \pm 1
\end{equation}
such that $(x , y, z)$ satisfies the system of quadratic ternary equations
\begin{equation}\label{twoquadraticsmodified}
\mathbf{Q}_{1}(X , Y, Z) = u , \, \, 
\mathbf{Q}_{2}(X , Y, Z) = v.
\end{equation}
\end{lemma}

\subsection{The trivial solution of $F(U, V) =1$}
First we notice that $F(U, V)$ is monic and therefore $(u , v) = (1, 0)$ satisfies the equation $F(U, V) = \pm 1$.
 This corresponds to the system of equations 
\begin{eqnarray}\label{system=0,1}\nonumber
\mathbf{Q}_1 (X, Y, Z)& =& 1\\ 
\mathbf{Q}_2 (X, Y, Z)& = & 0,
\end{eqnarray}
where the ternary quadratic forms $\mathbf{Q}_1$ and $\mathbf{Q}_2$ are defined in \eqref{defofQ1} and \eqref{defofQ2}.

A  special solution to the system of equations \eqref{system=0,1} is $(x , y, z) = (1, 0, 0)$ as $\xi$ is trivially a monogenizer of $\mathbb{Z}[\xi]$ (see \eqref{I0xyz}).

Assume  $x, y, z \in \mathbb{Z}$ satisfy \eqref{system=0,1}. Then
\begin{equation}\label{Q2=0}
\mathbf{Q}_{2}(x , y, z) = y^2 - xz - a_{1} yz + a_{2} z^2 = 0.
\end{equation}
If $z= 0$ then $y=0$. Since $x, y, z$ also satisfy $\mathbf{Q}_1(X , Y , Z) = 1$, we  conclude that $x=1$.

 Now  assume that $z\neq 0$. From \eqref{Q2=0}, we conclude that $z \mid y^2$. Let $q= \gcd(z, y)$, $y = q y'$ and $z =  q z' $, with $\gcd(y', z') =1$. We may rewrite \eqref{Q2=0} as 
\begin{equation*}
\mathbf{Q}_{2}(x , y, z) = y'^2q^2 - xz'q - a_{1} y'z' q^2 + a_{2} z'^2q^2 = 0
\end{equation*}
to conclude that $q \mid xz'$ and $z' \mid q$. Since $(x , y, z)$ satisfies the system \eqref{system=0,1}, in particular $\mathbf{Q}_1(x, y, z) =1$, we have $\gcd(q, x) = 1$ and therefore $q \mid z'$. So we have $z' = \pm q$ and   $q^2  = \pm z$. Since $(x , y, z)$ and $(-x, -y, -z)$ give the same monogenization, we may assume $z\geq 0$ and $q^2= z$. Now we can express $x, y$ and $z$ in terms of two integers $q$ and $p$ as follows:
\begin{equation}\label{GPPE8}
x = p^2  - a_{1} pq+ a_{2} q^2, \, \,  y= pq, \, \, z = q^2.
\end{equation}
The  parametrization \eqref{GPPE8} can be done for any $(x , y, z) \neq (1, 0, 0)$ that satisfies \eqref{Q2=0}. Substituting the parametrized values for variables $X, Y$ and $Z$ in \eqref{defofQ1}, we may express the ternary quadratic form $\mathbf{Q}_1(X, Y, Z)$ as a quartic binary form in variables $P$, $Q$,
where
\begin{equation}\label{XYZpq}
X(P, Q)= P^2  - a_{1} PQ+ a_{2} Q^2, \, Y(P, Q)= PQ, \, Z(P, Q)= Q^2.
\end{equation}
We note that each $X(P,Q)$, $Y(P,Q)$ and $Z(P,Q)$ is a binary quadratic form in variables $P$ and $Q$.  The parametrization \eqref{GPPE8} was considered for $z \neq 0$, however the trivial (and special) solution $(x, y, z) = (1,0,0)$ also corresponds to a solution of the quartic Thue equation
$$
\mathbf{Q}_{1}(X(P, Q) , Y(P,Q), Z(P,Q)) = 1,
$$
 namely $(p, q) = (1, 0)$.

Let us define the quartic binary form 
\begin{equation}\label{letusquarticThue}
\mathcal{Q}_{(1, 0)}(P, Q) =\mathcal{Q}(P, Q) = \mathbf{Q}_{1}(X(P, Q) , Y(P,Q), Z(P,Q)).
\end{equation}
We have shown that the number of solutions $(X, Y, Z) \in \mathbb{Z}^3$ of the system of ternary equations \eqref{system=0,1} is equal to the number of integer solutions $(p, q)$ of the quartic Thue equation 
$$
\mathcal{Q}(P, Q)  = 1.
$$

Via \eqref{system=0,1}, we may substitute  the parameter $X$ by  $\frac{Y^2 - a_1 YZ + a_2 Z^2}{Z^2}$  in $\mathbf{Q}_1(X, Y, Z)$ to get
$$
 \mathbf{Q}_1(X, Y, Z) = Z^4 \mathbf{P}\left(\frac{Y}{Z} - a_1\right),
 $$
 where 
 $\mathbf{P}(T)$ is the minimal polynomial of $\xi$ defined in \eqref{minpolydef}.  In other words,
 $$
\mathcal{Q}(P, Q) = Q^4 \mathbf{P}\left(\frac{P}{Q} - a_1\right).
$$
Since $a_1 \in \mathbb{Z}$, we conclude that  the discriminant of the quartic form $\mathcal{Q}(P, Q)$ is equal to the discriminant of $\xi$, and therefore, to the discriminant of the cubic form $F(U, V)$.

We also note that $
\mathcal{Q}(P, Q)$ is a monic binary form, i.e., the coefficient of the term $P^4$ equals $1$. This confirms the existence of the trivial solution $(p, q) = (1, 0)$ of the Thue equation
$\mathcal{Q}(P, Q) = 1$.

We conclude that the trivial solution $(1, 0)$ of the cubic Thue equation $F(U, V) =1$ corresponds to a quartic Thue equation,  namely 
$\mathcal{Q}_{(1, 0)}(P, Q) =1$, defined in \eqref{letusquarticThue}.
 Moreover, by \eqref{GPPE8}, each pair of solution $(p, q) \in \mathbb{Z}^2$  corresponds to the monogenizer
$$
X(p, q) \xi + Y(p, q) \xi^2 + Z(p, q) \xi^3
$$
of the order $\mathbb{Z}[\xi]$. Clearly, the monogenizer $\xi$ is produced by the solution
 $(p, q) = (1, 0)$ of the quartic Thue equation.

\subsection{Non-trivial solutions of $F(U, V) =1$} 
For non-trivial  solutions of the Thue equation  \eqref{cubicthue=1}, in \cite{GaPePo}  the system of ternary quadratic equations \eqref{Q'system}   is reduced to a quartic Thue equation with a parametrization similar to \eqref{XYZpq} (see  equations (8) and (10) of \cite{GaPePo}). We simplify such a parametrization with help of a $\textrm{GL}_2(\mathbb{Z})$ matrix that maps any given  primitive solution of a Thue equation to the trivial solution $(1 , 0)$ of an equivalent Thue equation.
More precisely, assume that $(u_0, v_0) \in \mathbb{Z}^2$, with $(u_0, v_0) \neq (1, 0)$, satisfies \eqref{cubicthue=1}. We have $\gcd(u_0, v_0) =1$ and therefore we may choose fixed  $s, t \in \mathbb{Z}$ so that 
$$
s u_0 + t v_0 = 1.
$$
Consequently, if $(x , y, z) \in \mathbb{Z}^3$ satisfies the system of equations in  \eqref{twoquadratics} with $(u, v) = (u_0, v_0)$,  then $(x, y, z)$ will satisfy
\begin{eqnarray}\label{Q'system}\nonumber
\mathbf{Q'}_1(X, Y, Z) & =& s\mathbf{Q}_1(X, Y, Z) + t \mathbf{Q}_2(X, Y, Z) = 1\\
\mathbf{Q'}_2(X , Y, Z) & =& v_0 \mathbf{Q}_1 - u_0 \mathbf{Q}_2 = 0.
\end{eqnarray}
The next step is to  express this system as an equation of a quartic binary form to $1$, via the parametrization  \eqref{XYZpq}.

Let 
$A = \left(\begin{array}{ll}
s & t\\
-v_0 & u_0
\end{array}\right) \in \textrm{GL}_2(\mathbb{Z})
$. Clearly we have 
$$A \left(\begin{array}{l}
u_0\\
v_0 \end{array} \right) = \left(\begin{array}{l}
1\\
0 \end{array} \right).
$$
The matrix $A^{-1} \in \textrm{GL}_2(\mathbb{Z})$ acts on the binary cubic form $F(U, V)$ to produce the equivalent binary form $F_{A^{-1}}(U, V)$.  The solution $(u_0, v_0)$ of $F(U, V) = 1$ corresponds  to the solution $(1 , 0)$ of the cubic equation $F_{A^{-1}}(U, V) = 1$.

Since $(1, 0)$ satisfies the equation $F_{A^{-1}}(U, V) = 1$, the cubic binary form $F_{A^{-1}}(U, V)$ is monic. Similar to \eqref{letusquarticThue}, and via parametrization \eqref{XYZpq}, we obtain the binary quartic form
\begin{equation}\label{letusquarticThueuv}
\mathcal{Q}_{(u_0, v_0)}(P, Q) = \mathbf{Q'}_{1}(X(P, Q) , Y(P,Q), Z(P,Q)),
\end{equation}
with $\mathbf{Q'}_1$  defined in \eqref{Q'system}. Therefore, in order to solve the system of ternary equations \eqref{Q'system} one can solve the quartic Thue equation 
\begin{equation}\label{Thueu0v0}
\mathcal{Q}_{(u_0, v_0)}(P, Q) = 1
\end{equation}
in integers $P, Q$.

\subsection{Conclusion}
Let $\xi$ be an algebraic integer of degree $4$ with the minimal polynomial given in \eqref{minpolydef}. In order to count the number of monogenizations of $\mathbb{Z}[\xi]$, we defined the integral cubic form $F(u, v)$ in \eqref{defofFu,v}, and the integral quadratic forms  $\mathbf{Q}_1(X, Y, Z)$ and $\mathbf{Q}_2(X, Y, Z)$ in \eqref{defofQ1} and \eqref{defofQ2}, respectively. The coefficients of these forms are all expressed in terms of the coefficients of the minimal polynomial of $\xi$.
We showed that the following three  numbers are equal:
\begin{enumerate}
\item the number of solutions to the cubic Thue equation $F(U, V) =\pm 1$ in  \eqref{cubicthue=1}, 
\item the number of systems of ternary quadratic equations \eqref{Q'system}, and
\item the number of quartic Thue equations \eqref{Thueu0v0}.
\end{enumerate}

We have also shown that for any fixed solution $(u, v) \in \mathbb{Z}^2$ of the  cubic Thue equation $F(U, V) =\pm 1$ in  \eqref{cubicthue=1}, each solution $(p, q)$ of the corresponding quartic Thue equation \eqref{Thueu0v0} provides a monogenizer $X(p, q) \xi + Y(p, q) \xi^2 + Z(p, q) \xi^3$, with the integral binary quadratic forms $X(P,Q)$, $Y(P, Q)$ and $Z(P, Q)$  defined in \eqref{XYZpq}.

Therefore, the number of monogenizations of $\mathbb{Z}[\xi]$ is bounded by an upper bound for the number of integer solutions to cubic Thue equations multiplied by an upper bound for the number of integer solutions  to quartic Thue equations. Proposition \ref{BenOka} provides upper bounds for the number of solutions of cubic Thue equations and Proposition \ref{main26}
provides upper bounds for the number of solutions of quartic Thue equations.

\section*{Acknowledgements} I am grateful to the anonymous referee for many helpful comments and suggestions, which improved an earlier version of this article significantly.

 I thank \emph{Professor Manjul Bhargava}
for inspiring and insightful conversations about the general topic of this article. I thank \emph{Professor K\'alm\'an Gy\H{o}ry} for encouragement, sharing his beautiful work, and comments  on an earlier version of this article.

During the completion of this project I visited the department of mathematics at Cornell University and was supported by Ruth I. Michler Memorial Prize. I am grateful to \emph{Professor Michler's family} and AWM for creating this invaluable opportunity and to Cornell math department for  hospitality. In particular, I  thank \emph{Professor Ravi Ramakrishna}  for warmly welcoming me and making me feel at home in Ithaca.

This research has been supported in part  by  \emph{the Simons Foundation Collaboration Grants}, Award Number 635880, and by  \emph{the National Science Foundation} Award DMS-2001281.


\end{document}